\documentclass{article}

\usepackage{graphicx} % Required for inserting images
\usepackage{float}
\usepackage{amssymb,amsbsy,amsthm,amsmath,epsfig,babel}
\usepackage{authblk}

\title{Structure and form of the solutions of the Erd\H{o}s-Straus conjecture}
\author[$\dagger$]{Miguel Angel Lopez}
\affil[$\dagger$]{Complutense University of Madrid,
Plaza de Ciencias 3,
Madrid 28040,
Spain}

\date{}

\begin{document}

\maketitle

\begin{center}email: migulo23@ucm.es, dagon.magnus@gmail.com\end{center}

\begin{abstract}
    Abstract: In this paper we classify certain values of p that satisfy the Erdos-Straus conjecture, concerning the decomposition of fractions of the form $\frac{4}{n}$ as sum of three fractions with numerator identically equal to 1, not according to their modular similarity but to the fact that they share solutions with identical structure. We classify all solutions that satisfy that they are of the form $(du,dv,duv)$  and find characterizations for values of $p$ that were initially difficult to classify, such as $p=1009.$  We also study the solutions $(x,y,z)$ that satisfy that $gcd(x,y,z)=x$ and, finally, we classify all the cases in which any of these two variables coincide with each other.
\end{abstract}

\textbf{Keywords:} Erdös-Straus, diophantine equation, unit fraction, egyptian fraction, congruences, integral solution, prime numbers

\textbf{Mathematics Subject Classification:} 11D72, 11A07, 11A41

\section{Introduction}

In ancient Egypt arithmetic relied on unit fractions, or fractions of numerator 1, to represent all other fractions. One of the uses of this kind of decomposition is to divide quantities in an intuitive and practical way. For example, if you want to divide 13 pizzas among 12 people, it is possible to give one pizza to each person and then break the last one into twelve pieces, which could wreck that last pizza in the process. However, since \begin{equation*} \frac{13}{12}=\frac{1}{2}+\frac{1}{3}+\frac{1}{4} \end{equation*} it will therefore be sufficient to split 6 pizzas into halves, 4 pizzas into thirds and 3 pizzas into quarters, and each person will eat exactly the same.
\\ \\
Over the centuries many discoveries have been made about the decomposition of fractions into sums of unit fractions, such as Fibonacci's Greedy Algorithm. In 1948, Paul Erdös and Ernst G. Straus formulated a conjecture that states the following: the equation \begin{equation} \frac{4}{n}=\frac{1}{x}+\frac{1}{y}+\frac{1}{z} \end{equation} has at least one solution where $x,$ $y,$ and $z$ are positive integers. There are many modular identities that solve this equation, for example for the case $n\equiv 2 \pmod3$ we can use the following expression: \begin{equation*} \frac{4}{n}=\frac{1}{n}+\frac{1}{(n+1)/3}+\frac{1}{n(n+1)/3} \end{equation*}

A simple but fundamental observation leads us to the fact that, in order to prove the conjecture, it is sufficient to do so for prime numbers, since every composite number is a multiple of a prime number and, therefore, if there exists a prime number p such that \begin{equation*} \frac{4}{p}=\frac{1}{x}+\frac{1}{y}+\frac{1}{z} \end{equation*} is a solution, then \begin{equation*} \frac{4}{kp}=\frac{1}{kx}+\frac{1}{ky}+\frac{1}{kz} \end{equation*} where $k$ is any positive integer. For this reason, and unless otherwise stated, we will henceforth restrict ourselves to the $p\in\mathbb{N}$ prime case.
\\ \\
 L. J. Mordell studied this equation among many others in \cite{mordell_diophantine_1970} and derived identities for the cases $p\equiv 3 \pmod4$, $2$ or $3\pmod5$, 3, 5 or $6\pmod7$ and $5\pmod8$. These combined identities solve all cases except those where $p$ is congruent to 1, 121, 169, 289, 361 or $529\pmod{840}$. Note that all these numbers are perfect squares, which is no coincidence. The first prime number for which the above relations fail to solve is 1009, a number we will discuss later in this paper.
\\ \\
Mordell himself wondered whether it might not be possible to find a sufficient number of identities such that all possible cases would be completely covered. This possibility was limited, however, by his discovery that if an identity exists for a set of values $p\equiv r \pmod{q}$ then $r$ cannot be a quadratic residue module $q$. For example, no such identity can exist for values of $p$ congruent to $1 \pmod{q}$  since 1 is always a quadratic residue module $q$ for any natural value of $q$. This implies, in fact, that it is not possible to find a value $q$ such that identities can be found for all elements of $\mathbb{Z}_q$.
\\ \\
The conjecture has been verified up to value $10^7$ by Yamamoto \cite{yamamoto_diophantine_1965} and $10^{14}$ by Swett \cite{salez_erdos-straus_2014}. Webb and others have shown that the natural density of possible counterexamples to the conjecture is zero, for as $N$ tends to infinity, the number of values in the interval $[1,N]$ that could be counterexamples tends to zero. Elsholtz and Tao have further shown that the number of solutions to the conjecture, averaged over the number of primes up to $n$ is bounded polylogarithmically in $n$ \cite{elsholtz_counting_2015}.
\\ \\
Recent authors who have worked on the conjecture are Zelator \cite{zelator_rational_2012}, Bradford \cite{bradford_note_2020} and Monks and Velingker \cite{monks_erdos-straus_nodate}, and I will rely on the work of some of them to continue and extend their ideas and theorems.

\section{Main body}

Before going into the theorems presented here, it is convenient to establish a notation. We will consider the values of each possible solution ordered in an increasing order, i.e., $x \leq y \leq z$. We will write a solution for each value $p\in N$ as $(x,y,z)$, what will allow us to list them in a more condensed form. We will use the following definition:
\\ \\
\textbf{Definition.} \textit{We say that a number $p\in N$ is Egyptian of order 3 if there exists a triplet of the form $(x,y,z)$ such that the Erdos-Straus conjecture is fulfilled for the fraction $\frac{4}{p}$.}
\\ \\
This notation will help us to establish some initial properties. Many of them are well known and even easy to deduce, but others are not so easy. I will refer to the work of Monks and Velingker \cite{monks_erdos-straus_nodate} and Zelator \cite{zelator_rational_2012}, where proof of these statements can be found, among many others that we will not use in this paper:
\\ \\
\textbf{Lemmas (Monks and Velingker).} \textit{Let p be prime and Egyptian of order 3 with solution $(x,y,z)$, then it holds:
\\ \\
(i) $p\nmid x$, $p \mid z$, $x<p$.
\\ \\
(ii) $x=y$ if and only if $p\equiv 3 \pmod4$ and $(x,y,z)=\left(\frac{p+1}{2},\frac{p+1}{2},\frac{p(p+1)}{4}\right)$.
\\ \\
(iii) $\lceil \frac{p}{4} \rceil \leq x \leq \frac{p+1}{2} $.}
\\ \\
We will further use this known result that, for example, Zelator proved in \cite{zelator_rational_2012}.
\\ \\
\textbf{Lemma (iv) (Zelator).} \textit{Let $p$ be prime, $n$ a positive integer, $n\geq 2$. We also assume that $gcd(n,p)=1$. We consider the diophantine equation in two variables}

\begin{equation*}
    nxy=p(x+y)
\end{equation*} 
\\ \\
\textit{Then this equation has a solution if and only if $n$ is a divisor of $p+1$. In fact, if $n\geq 3$, the only solution except order is equal to $(x,y)=\left(\frac{p+1}{n},p\left(\frac{p+1}{n}\right)\right)$.}
\\ \\
Note that the above equation is equivalent to the order 2 variant of the Erd\H{o}s-Straus Conjecture, that is,

\begin{equation*}
    \frac{n}{p}=\frac{1}{x}+\frac{1}{y}.
\end{equation*} 
\\ \\
The case $n=2$ has, in addition to the above solution, the trivial solution $(p,p)$.
\\ \\
Lemma ii) is particularly interesting, because although it is simple to prove, it does not seek to find an identity but to start from a property of the solution (that $x$ and $y$ are equal) and deduce which values can fulfill it. I will base my own statements on it.
\\ \\
An observation on the different solutions for each prime number $p$ reveals that there are many types of them but they can be classified into two large groups. Employing notation already used by Elsholtz and Tao first \cite{elsholtz_counting_2015}, and by Bradford later \cite{bradford_note_2020}, we will say that a solution $(x,y,z)$ is of Type I if $p\nmid x$, $p\nmid y$ and $p\mid z$, while a solution will be of Type II if $p\nmid x$, $p\mid y$ and $p\mid z$. If $p$ is not prime, there are other types of solutions.
\\ \\
The results of this paper are as follows.
\\ \\
\textbf{Theorem.} \textit{Let be $p \in \mathbb{N}$ prime, let $d \in \mathbb{N}$.
\\ \\
a) There exists a solution for $p$ with the form $(du,dv,duv)$ if and only if it exists $n \in \mathbb{N}$ such that $p\equiv -4d \pmod{4dn-1}$. In addition the values $u$ and $v$ are coprime, that is, $gcd(u,v)=1$.
\\ \\
b) Suppose that $p=4k+1$, there exists a solution for $p$ with the form $(du,dv,duv)$ if and only if it exists $n \in \mathbb{N}$ such that $k+1\equiv 3dn-d \pmod{4dn-1}$.
\\ \\
c) Suppose that $p=4k+1$, there exists a solution for $p$ with the form $(du,dv,duv)$ if and only if there exists $t\geq 0$ and a divisor $w$ of $k+1+t$ such that $w\equiv -1 \pmod{3+4t}$. In particular, if there exists a divisor $w$ of $k+1$ congruent to $2$ module $3$, $p$ is Egyptian of order $3$ with a solution of the form $(du,dv,duv)$.
\\ \\
d) There exists a solution $(x,y,z)$ for $p$ such that $gcd(x,y,z)=x$ if and only if there exists $n \in \mathbb{N}$, $n\geq\lceil\frac{p}{4}\rceil$ such that $4n-p\mid p+1$.
\\ \\
e) Suppose $p\geq3$. There exists a solution $(x,y,z)$ for $p$ such that $y=z$ if and only if $p\equiv 3 \pmod4$. Moreover, this solution is of the form $(u,2up,2up)$ with $u=\frac{p+1}{4}$.
\\ \\
f) Any prime number $p\geq3$, $p\equiv 3 \pmod4$, is Egyptian of order $3$ with exactly two solutions of the form $(x_1,x_1,z_1)$ and $(x_2,z_2,z_2)$. These two solutions are also related by means of the equalities $x_1=2x_2$, $2z_1=z_2$. If $p\not\equiv 3 \pmod4$, these two solutions never exist.}

\section{Notes on results}
Suppose $d=1$, then the condition of a) reduces to $p\equiv -4 \pmod{4n-1}$. This is an infinite succession of congruences beginning at $p\equiv 2 \pmod3$ and which in particular also encompasses other congruences not strictly contained in them, as for example $p\equiv 3 \pmod4$, since a solution always valid for $p\equiv -4\pmod{4n-1}$,  when $p$ is of the form $p=4k+3$, is to take $n=k+2$. These solutions are all of the form $(u,v,uv)$; e.g., $p=23$ has the solution $(8,23,184)$.
\\ \\
It is interesting to note that solutions of the form $(u,v,uv)$ also appear between values of $p$ of the form $p=4k+1$, as for example with $p=29$, having not only one but two of these solutions: $(10,29,290)$ and $(8,87,696)$.
\\ \\
If we extend the spectrum further and allow other values of $d$ different from $1$, the first surprise comes when the elusive value of $1009$, for which there was no congruence to solve it among those of Mordell, has two solutions of the form $(253,88792,2042216)$ and $(253,92828,1021108)$, with values of $d$ equal to $11$ and $23$, respectively. This implies that we have found two groups of congruences to which $1009$ belongs and which have a solution by means of an identity, and these are
\begin{equation*} p\equiv307\pmod{351}
\end{equation*}
\begin{equation*} p\equiv275\pmod{367}
\end{equation*}
These two congruences correspond to values of $dn$ equal to $88$ and $92$, and have been easily discovered by means of a simple Excel search using section a). It is worth noting that, according to the proof of this section, we will see that $p\mid v$, so all of these solutions fall into the Type II category.
\\ \\
The question arises whether every prime number $p$ always has a solution of the form $(du,dv,duv)$. The answer is negative: the first value not to have one is $p=193$. However, it is true that an immense number of values, even some that do not fall into other classes of congruences, do have solutions of this type. In particular, up to the value $p=4369$, only two prime numbers do not have a solution of this type: $193$ and $2521$.
\\ \\
Paragraph b) simplifies the accounts for the particular case we are most concerned with, that of $p\equiv 1\pmod4$.  Although there are many other more complex congruences, this one in particular is of great interest because of the very nature of the Erdös-Straus Conjecture, where the number $4$ is of special importance.
\\ \\
Paragraph c) is a result which, although it is an immediate consequence of paragraph b), offers fast and surprising calculation results even in the simple particular case of $t=0$.  For example, for the already mentioned case of $p=1009$, we have that $k=252$, and it is easy to verify that $k+1=21\cdot 23$.  Both divisors are congruent to 2 module 3 and, therefore, both serve to produce a solution of the type $(du,dv,duv)$, that are the two mentioned above. In the case of larger prime numbers the results can be even more immediate: for example, for the already considerable value of  $p=35617$, we have $k+1=8905$, that it has as obvious divisor 5 and, therefore, this proves automatically that it has a solution of type $(du,dv,duv)$ and even allows us to construct it following the steps of the demonstration, arriving at
\begin{equation*} (du,dv,duv)=(8095,126867754,634338770).\end{equation*}
The case $t=0$  covers a large number of values, but not all. The first value not to satisfy it is $p=73$,  where $k+1=19$,  that it has no divisors congruent to 2 module 3. The statement does hold, however, for the case $t=2$,  because $k+1+t=21$  and  $21\equiv -1 \pmod{3+4\cdot2}$.
\\ \\
It is straightforward to show that, if $p=4k+1$,  then a mandatory requirement for $p$ to not have solutions of the form $(du,dv,duv)$ is that $k\equiv 0 \pmod6$.  This occurs because if $k$ is odd, then $2\mid k+1$,  and therefore the cases $k\equiv1,3,5\pmod6$ are all covered by being all odd. The case $k\equiv2\pmod6$ does not apply because then $p$ would be divisible by 3 and therefore would not be prime, and in the case $k\equiv4\pmod6$  we have that, if $k=6u+4$,  then $k+1=6u+5$  and, therefore, is always congruent to 2 module 3 for every value of $u$. It can therefore be concluded that, if $p\not\equiv1\pmod{24}$, then $p$ always has at least one solution of the form  $(du,dv,duv)$.
\\ \\
In the particular case where $u=1$  we have $(x,y,z)=(d,dv,dv)$,  a solution type that satisfies two properties:

\begin{itemize}
    \item $gcd(x,y,z)=x$
    \item $y=z$
\end{itemize}

The first condition is characterized in section d), while the second is treated in sections e) and f). Section d) does not indicate the number of solutions $gcd(x,y,z)=x$ that $p$ can have. For example, for $p=29$, there are two such solutions, which are $(8,80,2320)$  and $(11,22,638)$  since there are two numbers $n\geq\lceil\frac{29}{4}\rceil=8$  that $4n-29\mid30$,  which are 8 itself and also 11.
\\ \\
In the case $p\equiv 3 \pmod4$ we have that the equation $4n-p=1$ is always solvable when $n=\frac{p+1}{4}$, so $p$ will always have at least two solutions such that $gcd(x,y,z)=x$. The proof concretizes its form, and says that they are $(n,2np,2np)$ and $(n,n(p+1),np(p+1))$,  with value of $n$ equal to that described above. Note that the first of these two solutions satisfies that $y=z$.  The following section focuses on that type of solution.
\\ \\
Regarding sections e) and f), for the case $p=2$ the answer as to whether there is a solution of the form $(x,y,z)$ with $y=z$  is trivial and positive with the solution $(1,2,2)$ which, moreover, is the only solution for this value. Section e) characterizes these solutions and shows that they are unique for each value of $p$ and appear in the case $p\equiv 3 \pmod4$.  Section f), finally, relates these solutions to those of the form $(x,y,z)$  with $x=y$,  showing that they always appear in pairs in the values  $p\equiv 3 \pmod4$.
\section{Proofs}
\textit{Proof of a).} First we will prove the implication from left to right. We have a solution that satisfies the equation
\begin{equation*}
    \frac{4}{p}=\frac{1}{du}+\frac{1}{dv}+\frac{1}{duv}
\end{equation*}
and which we can rewrite as
\begin{equation*}
    4duv=p(1+u+v)
\end{equation*}
We know from lemma i) that $p\mid duv$  and also $du<p$,  so $p\mid v$.  We denote $v=pv_1$,  introduced in the above equation we are left with
\begin{equation*}
    4dupv_1=p(1+u+pv_1),
\end{equation*}
and after simplifying we obtain that
\begin{equation*}
    4duv_1=1+u+pv_1.
\end{equation*}
By subtracting to obtain the value of $u$, we have that
\begin{equation*}
    u=\frac{1+pv_1}{4dv_1-1},
\end{equation*}
and this number will only be natural if $4dv_1-1\mid 1+pv_1$  or, what is the same thing, there exists a certain $n \in \mathbb{N}$  such that
\begin{equation*}
    1+np\equiv 0\pmod{4dn-1}.
\end{equation*}
Rearranging the congruence we have that
\begin{equation*}
    np\equiv -1\pmod{4dn-1}
\end{equation*}
and taking into account that $gcd(n,4dn-1)=1$  and that in fact $4dn\equiv 1 \pmod{4dn-1}$,  we multiply by $4d$ on both sides and we obtain that
\begin{equation*}
    p\equiv -4d\pmod{4dn-1},
\end{equation*}
which was what we were looking for.
\\ \\
To prove the converse implication, let $n \in \mathbb{N}$  be such that the above congruence is satisfied for a certain value $d \in \mathbb{N}$.  Again we know that $4d$ and $n$ are inverses of each other in $\mathbb{Z}_{4dn-1}$,  so we multiply by $n$ on both sides and obtain the equation
\begin{equation*}
    np\equiv -1\pmod{4dn-1}.
\end{equation*}
Therefore we obtain that
\begin{equation*}
    np+1\equiv 0\pmod{4dn-1}
\end{equation*}
and, therefore, we can conclude that the number $\frac{1+np}{4dn-1}\in \mathbb{N}$.  We call this number $u$ and define $v=np\in \mathbb{N}$.  Then, knowing that $u(4dn-1)=1+np$,  we arrive at the equation
\begin{equation*}
    4dnu=1+u+np,
\end{equation*}
where we first multiply by $p$ on both sides:
\begin{equation*}
    4dnup=p(1+u+np),
\end{equation*}
then we substitute  $v=np$:
\begin{equation*}
    4duv=p(1+u+v),
\end{equation*}
and finally we rearrange to obtain that
\begin{equation*}
    \frac{4}{p}=\frac{1}{du}+\frac{1}{dv}+\frac{1}{duv}.
\end{equation*}
To show that $u$ and $v$ are coprime, we define $m=gcd(u,v)$,  then $u=mu_1$, $v=mv_1$ and the above equation can thus be rewritten as
\begin{equation*}
    4dm^2u_1v_1=p(1+mu_1+mv_1).
\end{equation*}
A simple rearrangement leads to the following,
\begin{equation*}
    4dm^2u_1v_1-pmu_1-pmv_1=p,
\end{equation*}
from which it is immediate to deduce that $m\mid p$,  therefore $m=1$  or $m=p$.  If the latter occurs, then $mu_1\geq p$,  which implies that $dmu_1=du\geq p$,  but this is impossible because $du<p$.  Therefore it must occur that  $m=1$.
\\ \\
\textit{Proof of b).} We know from a) that a solution for $p$ with the form $(du,dv,duv)$  exists if and only if it exists $n \in \mathbb{N}$  such that $4k+1\equiv -4d\pmod{4dn-1}$.  Adding 3 to both sides we have that
\begin{equation*}
    4k+4\equiv 3-4d\pmod{4dn-1}
\end{equation*}
As we know that the inverse of 4 in $\mathbb{Z}_{4dn-1}$  is $dn$, we multiply on both sides and we have that
\begin{equation*}
    (4dn)k+4dn \equiv 3dn-(4dn)d\pmod{4dn-1}
\end{equation*}
This is equivalent to
\begin{equation*}
    k+1 \equiv 3dn-d\pmod{4dn-1}
\end{equation*}
which is exactly what we wanted to demonstrate.
\\ \\
\textit{Proof of c).} We define $d=\frac{k+1+t}{w}$.  We know that we can write $w=(3+4t)w_1-1$,  with $w_1$  natural and strictly positive. We define $n=w_1$.  It is clear then that
\begin{equation*}
    k+1+t=wd=((3+4t)n-1)d=(3+4t)dn-d
\end{equation*}
\begin{equation*}
    k+1=(3+4t)dn-d-t=3dn-d+t(4dn-1),
\end{equation*}
and finally, by module $4dn-1$  we arrive at
\begin{equation*}
    k+1\equiv 3dn-d \pmod{4dn-1}
\end{equation*}
Therefore we fall under the hypothesis of b) and $p$ has a solution of type $(du,dv,duv)$. For the opposite implication, it suffices to see that all the above steps are reversible.
\\ \\
\textit{Proof of d).} Let's first look at the left implication. We can write the solution as $(a,ab,ac)$,  where $a,b,c$  are natural numbers. Then
\begin{equation*}
    \frac{4}{p}=\frac{1}{a}+\frac{1}{ab}+\frac{1}{ac},
\end{equation*}
Moving $\frac{1}{a}$ to the left side and multiplying both sides by $a$, we have that
\begin{equation*}
    \frac{4a-p}{p}=\frac{1}{b}+\frac{1}{c}.
\end{equation*}
Thanks to lemma iii) we know that $a\geq \lceil\frac{p}{4}\rceil$ and therefore $4a-p\geq 1$.  We choose then $n=a$.  If $4n-p=1$,  then $4n-p\mid p+1$  trivially. If $4n-p\geq 2$,  using lemma iv), we have that $4n-p\mid p+1$.
\\ \\
Suppose now that we have a value  $n \in \mathbb{N}$, $n\geq\lceil\frac{p}{4}\rceil$,  such that $4n-p\mid p+1$.  This implies that $4n-p\geq 1$.  We distinguish two cases:
\\ \\
Case 1. $4n-p=1$.  This situation is only possible when $p$ is of the form $4n-1$,  i.e., if $p\equiv 3 \pmod4$.  In this case, it is a classical result that the fraction $\frac{1}{p}$  can be written as a sum of two unit fractions in two ways, which are
\begin{equation*}
    \frac{1}{p}=\frac{1}{2p}+\frac{1}{2p}=\frac{1}{p+1}+\frac{1}{p(p+1)}
\end{equation*}
and therefore, as $4n-p=1$,  we have that
\begin{equation*}
    \frac{4n-p}{p}=\frac{4n}{p}-1=\frac{1}{2p}+\frac{1}{2p}=\frac{1}{p+1}+\frac{1}{p(p+1)}
\end{equation*}
Performing elementary arithmetic, we obtain that $p$ has two possible solutions, $(n,2np,2np)$  and  $(n,n(p+1),np(p+1)$.
\\ \\
Case 2. $4n-p\geq 2$.  We know that $4n-p\mid p+1$,  therefore we apply lemma iv) and we have that
\begin{equation*}
    \frac{4n-p}{p}=\frac{1}{b}+\frac{1}{c}
\end{equation*}
where in fact $(b,c)=\left(\frac{p+1}{4n-p},p\left(\frac{p+1}{4n-p}\right)\right)$ and there also exists the additional solution $(b,c)=(p,p)$  for the case $4n-p=2$.  Again using elementary arithmetic, we find that the solution is of the form
\begin{equation*}
    (x,y,z)=\left(n,n\left(\frac{p+1}{4n-p}\right),np\left(\frac{p+1}{4n-p}\right)\right)
\end{equation*}
if $4n-p\geq 2$,  and also
\begin{equation*}
    (x,y,z)=(n,np,np)
\end{equation*}
if $4n-p=2$.  Since all possible cases are covered, there is always a solution of this form.
\\ \\
\textit{Proof of e).} If $p\equiv 3 \pmod4$, it suffices to use the solution proposed in the statement of e) and perform a simple calculation to show that this solution is valid. Let us assume therefore that $p\equiv 3 \pmod4$  and we will check that it is the only possible solution.
\\ \\
Suppose that there exists a solution for $p$ with the requested form, then it must fulfill that
\begin{equation*}
    \frac{4}{p}=\frac{1}{x}+\frac{1}{y}+\frac{1}{y}
\end{equation*}
and therefore this implies that $4xy=p(2x+y)$.  By lemma i) we know that $p\mid z=y$,  we call then $y=py_1$.  We introduce it to obtain that $4xpy_1=p(2x+py_1)$  and, therefore, $4xy_1=2x+py_1$.  This instantly implies that $x\mid py_1$.  If $x=1$  then $x\mid y_1$  trivially. If $x>1$  we know that $gcd(x,p)=1$  and, therefore, $x$ does not divide $p$ and must divide $y_1$,  thanks to Euclid's Lemma. We call therefore $y_1=xy_2$.  Transferring this information to the equation we have that $4x^2y_2=2x+pxy_2$,  which implies that $4xy_2=2+py_2$.  Rearranging and taking out the common factor, we arrive at that
\begin{equation*}
    y_2(4x-p)=2
\end{equation*}
We know that $y_2 \in \mathbb{N}$, $4x-p \in \mathbb{N}$  by lemma iii). This automatically implies that $4x-p=1$  or  $4x-p=2$. We distinguish between the two cases.
\\ \\
Case 1. $4x-p=2$.  This implies that $p\equiv 2 \pmod4$  and we have already ruled out this case in the statement of the theorem, although it is worth noting that, if we were to develop it, we would arrive at the aforementioned solution $(1,2,2)$  applicable only to the value of  $p=2$.
\\ \\
Case 2. $4x-p=1$.  This implies that $p\equiv 3 \pmod4$.  In addition, $y_2=2$,  which implies that $y_1=2x$  and, finally, $y=2xp$.  Therefore our solution is of the form $(u,2up,2up)$,  where
\begin{equation*}
    \frac{4}{p}=\frac{1}{u}+\frac{1}{2up}+\frac{1}{2up}=\frac{2p+2}{2up}=\frac{p+1}{up},
\end{equation*}
which implies that
\begin{equation*}
    4=\frac{p+1}{u}
\end{equation*}
and finally
\begin{equation*}
    u=\frac{p+1}{4}.
\end{equation*}
\textit{Proof of f).} The uniqueness of the solution with $y=z$  is deduced by construction in e), and the uniqueness of the solution with $x=y$  was proved by Monks and Velingker in \cite{monks_erdos-straus_nodate} and is lemma iii). By construction, these two solutions have the respective forms
\begin{equation*}
    (x_1,y_1,z_1)=\left(\frac{p+1}{2},\frac{p+1}{2},\frac{p(p+1)}{4}\right)
\end{equation*}
\begin{equation*}
    (x_2,y_2,z_2)=\left(\frac{p+1}{4},\frac{p(p+1)}{2},\frac{p(p+1)}{2}\right)
\end{equation*}
A simple inspection leads us to $x_1=2x_2$, $2z_1=z_2$,    what we wanted to check. It has already been shown in e) for the second type, and by Monks and Velingker for the first type, that these two solutions can only exist when  $p\equiv 3 \pmod4$.

\section{Conclusion}

We have analyzed, for prime $p$, particular solutions to the Erdös-Straus Conjecture focusing on the form of their solutions. This has led us to systems of congruences or to new properties related to finding divisors of certain numbers. We have also proposed a simple theorem that demonstrates the potential of this approach by trivially solving numbers that until now were complex to study, such as 1009, the first value that Mordell failed to classify. These results cut across the classical approach, for instead of assuming a modular type of relation to draw a conclusion from it, they look for structural properties of the solutions to group the numbers according to that common structure.

\bibliographystyle{unsrt}
\bibliography{structure}
\end{document}